\documentclass[12pt,a4paper]{article}
\usepackage[T1]{fontenc}
\usepackage[english]{babel}
\usepackage{lmodern}
\usepackage{amssymb,amsmath,amsthm,mathtools}
\usepackage[hmargin=2.5cm,vmargin=3cm]{geometry}
\usepackage{microtype}
\usepackage{authblk}
\usepackage{cite}
\usepackage{graphicx,xcolor}
\usepackage{enumerate,enumitem}
\usepackage{hyperref}
\usepackage[capitalise,nameinlink,noabbrev]{cleveref}

\newtheorem{theorem}{Theorem}
\newtheorem{lemma}{Lemma}
\newtheorem{proposition}{Proposition}
\newtheorem{corollary}{Corollary}
\newtheorem{conjecture}{Conjecture}

\Crefname{conjecture}{Conjecture}{Conjectures}

\newcommand{\mad}{\operatorname{mad}}
\newcommand{\Poc}{P_{\mathrm{oc}}}
\newcommand{\R}{\mathbb{R}}

\begin{document}
\title{Weighted odd-cycle transversals and fractional colorings of planar graphs}

\author[1]{Xiaolan Hu}
\author[2]{Jiaao Li}
\author[2]{Zhouningxin Wang}

\affil[1]{\small School of Mathematics and Statistics and Hubei Key Laboratory of Mathematical Sciences, Central China Normal University, Wuhan 430079, China.}
\affil[2]{\small School of Mathematical Sciences and LPMC, Nankai University, Tianjin 300071, China.\linebreak
Emails: xlhu@ccnu.edu.cn; \{lijiaao, wangzhou\}@nankai.edu.cn.}

\date{}
\maketitle

\begin{abstract}
We bound fractional chromatic numbers by establishing a connection with low-weight odd-cycle transversals. Let $\ell\geq 3$ be an integer and let $G$ be a graph of odd girth at least $\ell$ whose odd-cycle covering polyhedron
is integral. For every nonnegative weight vector $\boldsymbol{w}=(w_v)_{v\in V(G)}\in\R_{\geq0}^{V(G)}$, we obtain a set $Z\subseteq V(G)$ such that $G-Z$ is bipartite and $w(Z)\leq\frac{\mad(G)}{2\ell}w(V(G))$. Here, $w(S)=\sum_{v\in S}w_v$ for any $S\subseteq V(G)$, and $\mad(G)$ denotes the maximum average degree of $G$. Consequently, if $\mad(G)<2\ell$, then $G$ has fractional chromatic number at most $\frac{4\ell}{2\ell-\mad(G)}$.

As applications, we prove that, for every positive integer $k$, every $n$-vertex planar graph of girth at least $2k+3$ has fractional chromatic number at most $\frac{2k+1}{k}(1-\frac{1}{kn+1})$, establishing the girth version of Naserasr's fractional coloring conjecture with a strict inequality. For $k=1$, we furthermore prove that the upper bound $\frac{3n}{n+1}$ holds for every $n$-vertex triangle-free planar graph and this is attained for infinitely many values of $n$.
\end{abstract}



\section{Introduction}\label{sec:intro}

All graphs considered are finite and simple. The \emph{girth} of a graph is the length of a shortest cycle, and its \emph{odd girth} is the length of a shortest odd cycle. For a positive integer $a$, let $[a]=\{1,\ldots,a\}$. For positive integers $a$ and $b$ with $a\ge b$, an \emph{$(a\colon b)$-coloring} of a graph $G$ is a mapping $\varphi:V(G)\to\binom{[a]}b$ assigning disjoint sets to adjacent vertices. The \emph{fractional chromatic number} of a graph $G$, denoted $\chi_f(G)$, is defined as follows:
$$
    \chi_f(G):=\inf\left\{\frac{a}{b}:
    G\text{ admits an }(a\colon b)\text{-coloring}\right\}.
$$
See Stahl~\cite{Stahl} and
Scheinerman--Ullman~\cite{ScheinermanUllman} for background.

For planar graphs, Naserasr~\cite{Naserasr} proposed the following
fractional coloring conjecture.

\begin{conjecture}{\em \cite{Naserasr}}
\label{conj:Naserasr-fractional}
For every integer $k\ge1$, every planar graph $G$ of odd girth at least $2k+3$ satisfies $\chi_f(G)\leq \frac{2k+1}{k}$.
\end{conjecture}

Naserasr~\cite{Naserasr} also noted that the condition cannot be weakened to that of
odd girth at least $2k+1$: for every positive integer $k$, there exists a planar graph $G$ of odd girth $2k+1$ with
$\chi_f(G)>\frac{2k+1}{k}$.
A stronger coloring assertion, contained in~\cite[Conjecture~3.3]{Naserasr}, is that every planar graph of odd girth
at least $2k+3$ admits a $(2k+1\colon k)$-coloring.

The stronger coloring assertion for $k=1$ follows from Gr\"otzsch's $3$-color theorem for triangle-free planar graphs. For $k=2$, Dvo\v{r}\'ak, \v{S}krekovski, and
Valla~\cite{DvorakSkrekovskiValla} proved that odd girth at least $9$
suffices for a $(5\colon2)$-coloring. For general $k$, Pirnazar and Ullman~\cite{PirnazarUllman} proved the
$(2k+1\colon k)$-colorability of planar graphs of girth at least $8k-4$,
while subsequent results of Lov\'asz, Thomassen, Wu, and
Zhang~\cite{LovaszThomassenWuZhang} yielded the same conclusion for planar graphs of odd girth at least $6k+1$. For the class of $K_4$-minor-free graphs, Beaudou, Foucaud, and Naserasr~\cite{BeaudouFoucaudNaserasr} proved that odd girth at least $2k+1$ suffices for a $(2k+1\colon k)$-coloring. Hu and Li~\cite{HuLi} obtained bounds on the circular chromatic numbers and fractional chromatic numbers of planar graphs with no cycles of lengths in specified intervals.

A related line of work studies fractional colorings under density conditions without assuming planarity. The \emph{maximum average degree} of $G$ is $\mad(G)=\max\{ \frac{2|E(H)|}{|V(H)|}: H\subseteq G\}.$
The Chen--Raspaud conjecture~\cite{ChenRaspaud} asserts that, for every integer $k\geq 2$, every graph $G$ of odd girth at least $2k+1$ with $\mad(G)<\frac{2k+1}{k}$ admits a $(2k+1\colon k)$-coloring. Chen and Raspaud~\cite{ChenRaspaud} proved the case $k=2$, and \L{}yczek, Nazarczuk, and Rz\k{a}\.{z}ewski~\cite{LyczekNazarczukRzazewski} established the case $k=3$.
Choi~\cite{Choi} proved the case $k=4$ and obtained simpler proofs for $k\in \{2,3\}$. Wu and Lu~\cite{WuLu} subsequently announced a proof of the conjecture, yielding $(2k+1\colon k)$-colorability for planar graphs of girth at least $4k+2$.

In this paper, we prove the girth version of \Cref{conj:Naserasr-fractional} with a strict inequality and obtain a quantitative bound in terms of the number of vertices. 
For $n$-vertex triangle-free planar graphs, we further establish the sharp fractional chromatic bound $\frac{3n}{n+1}$.
Our proof uses a weighted bound on odd-cycle transversals in terms of odd girth and maximum average degree, which applies more generally to graphs whose odd-cycle covering polyhedron is integral.

To state our main results, we need the following terminology. An \emph{odd-cycle edge cover} of $G$ is an edge subset $B\subseteq E(G)$ such that $G-B$ is bipartite, while an \emph{odd-cycle
transversal} of $G$ is a vertex subset $Z\subseteq V(G)$ such that $G-Z$ is bipartite.
Let $\R_{\geq 0}^{V(G)}$ denote the set of vectors $\boldsymbol{w}=(w_v)_{v\in V(G)}$ with all $w_v$ nonnegative. For any $\boldsymbol{w}\in\R_{\geq 0}^{V(G)}$
and $S\subseteq V(G)$, let $w(S):=\sum_{v\in S}w_v$. The \emph{odd-cycle covering polyhedron} of a graph $G$ is
$$
    \Poc(G):=\Big\{\boldsymbol{x}\in\R_{\ge0}^{E(G)}:
    \sum_{e\in E(C)}x_e\ge1\text{ for every odd cycle }C\text{ of }G\Big\}.
$$
An \emph{extreme point} of $\Poc(G)$ is a vector that cannot be written as $\lambda\boldsymbol{u}+(1-\lambda)\boldsymbol{v}$ for distinct $\boldsymbol{u},\boldsymbol{v}\in\Poc(G)$ and $0<\lambda<1$. The polyhedron $\Poc(G)$ is \emph{integral} if each of its extreme points has integer coordinates.

By the linear programming formulation of fractional chromatic numbers
(see \Cref{sec:cover} and \Cref{lem:criterion}), to obtain an upper bound on $\chi_f(G)$, it suffices   to find an independent
set of sufficiently large weight for every vector
$\boldsymbol{w}\in\R_{\geq0}^{V(G)}$.
Deleting a low-weight odd-cycle transversal leaves a bipartite graph, whose heavier bipartition part is an independent set containing at least half of the remaining total weight.
To find such a transversal, we assign each edge $uv$ the cost $\min\{w_u,w_v\}$ and select an endpoint of minimum weight from each edge of a low-cost odd-cycle edge cover. Maximum average degree controls the total edge cost, giving the
following bound on the weight of the transversal.

\begin{proposition}\label{thm:oct}
Let $\ell$ be an integer with $\ell\geq 3$ and let $G$ be a graph of odd girth at least $\ell$ such that $\Poc(G)$ is integral. For every $\boldsymbol{w}\in\R_{\geq 0}^{V(G)}$, there exists an odd-cycle transversal $Z\subseteq V(G)$
such that $$w(Z)\leq\frac{\mad(G)}{2\ell}w(V(G)).$$
\end{proposition}

Applying~\Cref{thm:oct} with the weighted characterization yields the following upper bound on the fractional chromatic number in terms of maximum average degree and odd girth.

\begin{theorem}\label{thm:general}
Let $\ell$ be an integer with $\ell\geq 3$. Every graph $G$ of odd girth at least $\ell$ with $\mad(G)<2\ell$ and integral $\Poc(G)$ satisfies $$\chi_f(G)\leq \frac{4\ell}{2\ell-\mad(G)}.$$
\end{theorem}

Both \Cref{thm:oct} and \Cref{thm:general} are sharp for odd cycles; see~\Cref{sec:cover} and~\Cref{sec:proofs}.

Guenin's characterization~\cite{Guenin} of weakly bipartite signed graphs, for which Schrijver~\cite{Schrijver} gave a shorter proof, implies that $\Poc(G)$ is integral whenever $G$ is $K_5$-minor-free.
We explain this consequence in \Cref{sec:cover}.
In particular, $\Poc(G)$ is integral for every planar graph $G$.
Combining~\Cref{thm:general} with the standard Euler formula estimate for planar graphs of given girth (see, \cite{Diestel2005}, for example) yields the following quantitative form of the girth version of~\Cref{conj:Naserasr-fractional}.

\begin{theorem}\label{thm:planar}
For all positive integers $k$ and $n$, every $n$-vertex planar graph $G$ of girth at least $2k+3$ satisfies $\chi_f(G)\leq \frac{(2k+1)n}{kn+1}$.
\end{theorem}

Since $\frac{(2k+1)n}{kn+1}=\frac{2k+1}{k}(1-\frac{1}{kn+1})<\frac{2k+1}{k}$, \Cref{thm:planar} establishes the desired strict inequality.
For $k=1$, it gives the bound $\frac{3n}{n+1}$ for planar graphs of girth at least $5$.
We next extend this bound to all triangle-free planar graphs.

Dvo\v{r}\'ak, Sereni, and Volec~\cite{DvorakSereniVolecPlanar} proved that every $n$-vertex triangle-free planar graph $G$ satisfies $\chi_f(G)\leq\frac{9n}{3n+1}$. They also established the stronger bound $\frac{3n}{n+1}$ when the maximum degree is at most four. Subsequently, Dvo\v{r}\'ak and Hu~\cite{DvorakHu} improved the general bound to $\frac{6n}{2n+1}$ by proving that every $n$-vertex triangle-free planar graph admits a $(6n\colon2n+1)$-coloring. Our next result improves both general bounds and extends the bound $\frac{3n}{n+1}$ to all triangle-free planar graphs.

\begin{theorem}\label{thm:triangle-planar}
Every $n$-vertex triangle-free planar graph $G$ satisfies $\chi_f(G)\leq \frac{3n}{n+1}.$ Moreover, the bound is tight.
\end{theorem}

The proof combines \Cref{thm:planar} with folding along facial $4$-cycles and a gluing argument along separating $4$-cycles. Extendable $3$-colorings of the boundary are used to match the
color-set intersections on the two sides.

We remark that the bound in \Cref{thm:triangle-planar} is sharp for infinitely many values of $n$. Indeed, for every integer $n\geq2$ with $n\equiv2\pmod3$, Jones~\cite{Jones} constructed an $n$-vertex triangle-free planar graph $G_n$ of maximum degree at most $4$ such that $\alpha(G_n)=\frac{n+1}{3}$, where $\alpha(G_n)$ denotes the independence number; see also the discussion preceding Theorem~2 of~\cite{DvorakSereniVolecPlanar} and the introduction of~\cite{DvorakHu}. Since $\chi_f(H)\geq\frac{|V(H)|}{\alpha(H)}$ for every graph $H$, we have $\chi_f(G_n)\geq\frac{n}{\alpha(G_n)}=\frac{3n}{n+1}.$ Together with~\Cref{thm:triangle-planar}, this yields $\chi_f(G_n)=\frac{3n}{n+1}$.

\medskip
The same weighted bound also gives fractional coloring results for $K_5$-minor-free graphs under odd-girth and density conditions. In particular, \Cref{thm:general} yields the following consequence.

\begin{corollary}\label{thm:K5minorfree}
For every positive integer $k$, every $K_5$-minor-free graph $G$ of odd girth at least $2k+3$ and $\mad(G)\leq \frac{2(2k+3)}{2k+1}$ satisfies
$\chi_f(G)\leq \frac{2k+1}{k}$.
\end{corollary}

Applying some known edge-density bounds for $K_5$-minor-free graphs, we obtain the following corollaries.

\begin{corollary}\label{cor:minor-bounds}
Let $\ell$ and $n$ be two integers with $\ell\geq 5$ and $n\geq 4$, and let $G$ be an $n$-vertex $K_5$-minor-free graph of odd girth at least $\ell$. The following statements hold.
\begin{enumerate}[label=(\arabic*)]
\item\label{cor:minor-odd-girth}
$\chi_f(G)\leq \frac{2\ell n}{(\ell-3)n+6}.$
\item\label{cor:minor-no-four}
If $G$ has no $4$-cycles, then $\chi_f(G)\leq \frac{10\ell n}{(5\ell-9)n+20}.$
\end{enumerate}
\end{corollary}

The remainder of the paper is organized as follows. \Cref{sec:cover} presents the weighted independent-set criterion and the edge-cover estimate, and proves~\Cref{thm:oct}.
\Cref{sec:proofs} proves~\Cref{thm:general} and its consequences and, in particular, provides a proof of the sharp bound for triangle-free planar graphs. \Cref{sec:conclusion} discusses further directions toward~\Cref{conj:Naserasr-fractional} and some open problems.

\section{Odd-cycle transversals and fractional colorings}\label{sec:prelim}\label{sec:cover}

Let $\R_{\geq 0}^{E(G)}$ denote the set of vectors $\boldsymbol{x}=(x_e)_{e\in E(G)}$ with all $x_e$ nonnegative. 
For any $\boldsymbol{x}\in\R_{\geq 0}^{E(G)}$ and $F\subseteq E(G)$, write
$x(F):=\sum_{e\in F}x_e$. 
For any $\boldsymbol{x}, \boldsymbol{y}\in\R_{\geq 0}^{E(G)}$, let $\boldsymbol{x}\cdot\boldsymbol{y}:=\sum_{e\in E(G)}x_ey_e$
denote the standard inner product in $\R_{\geq 0}^{E(G)}$.

Let $\mathcal{I}(G)$ denote the family of all independent sets of $G$. A \emph{fractional coloring} of $G$ is a vector $\boldsymbol{\lambda}=(\lambda_I)_{I\in\mathcal{I}(G)}\in\R_{\geq 0}^{\mathcal{I}(G)}$ such that
$\sum_{I\ni v}\lambda_I\geq1$ for every $v\in V(G)$. Its total weight is $\sum_{I\in\mathcal{I}(G)}\lambda_I$, and the fractional chromatic number is the optimum of the following primal and dual linear programs (for example, see~\cite{ScheinermanUllman}):
\par\noindent
\begin{minipage}[t]{0.48\linewidth}
\vspace{0pt}
\small
\begin{equation*}
\begin{aligned}
\chi_f(G)=\min\quad
 &\sum_{I\in\mathcal{I}(G)}\lambda_I\\
\text{subject to}\quad
 &\sum_{\substack{I\in\mathcal{I}(G)\\v\in I}}\lambda_I\geq1
 &&(v\in V(G)),\\
 &\lambda_I\geq0
 &&(I\in\mathcal{I}(G)).
\end{aligned}
\end{equation*}
\end{minipage}\hfill
\begin{minipage}[t]{0.48\linewidth}
\vspace{0pt}
\small
\begin{equation*}
\begin{aligned}
\chi_f(G)=\max\quad
 &w(V(G))\\
\text{subject to}\quad
 &w(I)\leq1
 &&(I\in\mathcal{I}(G)),\\
 &w_v\geq0
 &&(v\in V(G)).
\end{aligned}
\end{equation*}
\end{minipage}
\par

\bigskip
The following weighted characterization follows from~Theorem~1.4.1~of~\cite{ScheinermanUllman}, applied to the hypergraph on $V(G)$ whose hyperedges are independent sets of $G$; see also~\cite{DvorakSereniVolec}.

\begin{proposition}{\em \cite{ScheinermanUllman,DvorakSereniVolec}}
\label{lem:criterion}
Let $r$ be a positive rational number. A graph $G$ satisfies $\chi_f(G)\le r$ if and only if, for every vector $\boldsymbol{w}\in\R_{\geq 0}^{V(G)}$, there exists an
independent set $I\subseteq V(G)$ such that $w(I)\ge \frac{1}{r}w(V(G))$.
\end{proposition}

Applying Guenin's theorem~\cite{Guenin} (also by Schrijver~\cite{Schrijver}) to the signed graph obtained by assigning every edge of $G$ negative yields the following result, since a cycle contains an odd number of negative edges if and only if it is of odd length.
\begin{proposition}{\em \cite{Guenin,Schrijver}}\label{lem:integral}
If $G$ is $K_5$-minor-free, then $\Poc(G)$ is integral.
\end{proposition}

The following lemma provides the edge-cover estimate used in the proof of~\Cref{thm:oct}.

\begin{lemma}\label{lem:cover}
Let $\ell$ be an integer with $\ell\geq 3$ and let $G$ be a graph of odd girth at least $\ell$ such that $\Poc(G)$ is integral. For every vector $\boldsymbol{c}\in\R_{\geq 0}^{E(G)}$, there is an odd-cycle edge cover $B\subseteq E(G)$ such that $$c(B)\le\frac{1}{\ell}\,c(E(G)).$$
\end{lemma}

\begin{proof}
Let $\boldsymbol{c}\in\R_{\geq 0}^{E(G)}$ be a vector. We shall minimize $\boldsymbol{c}\cdot\boldsymbol{x}=\sum_{e\in E(G)}c_ex_e$ over
all $\boldsymbol{x}\in\Poc(G)$.
The polyhedron is nonempty and contained in the nonnegative
orthant, so it is pointed (that is, it contains no complete straight line). A standard result in linear programming
(see Theorems~2.6 and~2.8 of~\cite{BertsimasTsitsiklis})
states that a linear objective function bounded below on
a nonempty pointed polyhedron attains its minimum at
an extreme point.

Since the objective is bounded below by zero, there is an optimal extreme point. Since $\Poc(G)$ is integral, every extreme point of $\Poc(G)$ has integer coordinates. Let $\boldsymbol{x}^*$ denote the constant vector with $x_e:=\frac{1}{\ell}$ for every $e\in E(G)$. Since every odd cycle of $G$ has at least $\ell$ edges, $\boldsymbol{x}^*\in \Poc(G)$. By integrality, such an optimum $\boldsymbol{x}$ has nonnegative integer coordinates and satisfies $$\boldsymbol{c}\cdot\boldsymbol{x}\leq\boldsymbol{c}\cdot\boldsymbol{x}^*=\frac{1}{\ell}c(E(G)).$$

Let $B=\{e\in E(G):x_e\geq 1\}$.
The constraint for each odd cycle forces that cycle to meet $B$, so $B$ is an odd-cycle edge cover. Indeed, the coordinates on an odd cycle sum to at least $1$; since they are nonnegative integers, at least one is at least $1$. Consequently, deleting $B$ destroys every odd cycle, leaving a bipartite graph.
Since $\boldsymbol{c}$ has nonnegative coordinates,
$c(B)\leq\boldsymbol{c}\cdot\boldsymbol{x} \leq\frac{1}{\ell}c(E(G))$.
\end{proof}

Now we are ready to present the proof of~\Cref{thm:oct}.

\begin{proof}[Proof of~\Cref*{thm:oct}]
Let $G$ be a graph of odd girth at least $\ell$ such that $\Poc(G)$ is integral. 
For any vector $\boldsymbol{w}\in\R_{\geq 0}^{V(G)}$, we define a vector $\boldsymbol{c}\in\R_{\geq 0}^{E(G)}$ with $c_{uv}:=\min\{w_u,w_v\}$. By~\Cref{lem:cover}, there is an odd-cycle edge cover $B$ such that $c(B)\leq \frac{1}{\ell}c(E(G))$. We select an endpoint of minimum weight from each edge of $B$, and let $Z$ be the set of all selected vertices. It is easy to see that $G-Z$ is a subgraph of the bipartite graph $G-B$. Thus $Z$ is an odd-cycle transversal. By construction, $w(Z)\le c(B)\le\frac{1}{\ell}\sum_{uv\in E(G)}\min\{w_u,w_v\}.$

Moreover, let $b_1<\cdots<b_r$ list the distinct positive coordinates of $\boldsymbol{w}$, let $b_0=0$, and define $U_i=\{v\in V(G):w_v\ge b_i\}$ for $1\le i\le r$. A vertex of weight $b_j$ belongs to precisely $U_1,\ldots,$ and $U_j$, so $ \sum_{i=1}^r(b_i-b_{i-1})|U_i|=w(V(G)).$
Also, an edge $uv$ belongs to $E(G[U_i])$ if and only if $\min\{w_u,w_v\}\ge b_i$. Therefore,
{\small $$\sum_{uv\in E(G)}\min\{w_u,w_v\}=\sum_{i=1}^r(b_i-b_{i-1})|E(G[U_i])|
\leq\frac{\mad(G)}{2}\sum_{i=1}^r(b_i-b_{i-1})|U_i|
=\frac{\mad(G)}{2}w(V(G)).$$}%
Combining these inequalities, we obtain $w(Z)\leq\frac{\mad(G)}{2\ell}w(V(G))$.
\end{proof}

Note that the bound in~\Cref{thm:oct} is attained by odd cycles. Indeed, let $\ell$ be an odd integer with $\ell\geq 3$ and let $C_\ell$ be a cycle of length $\ell$. We assign weight $1$ to every vertex. Every odd-cycle transversal has weight at least $1$, while $\frac{\mad(C_\ell)}{2\ell}w(V(C_\ell))=1$.

\section[Fractional chromatic number]{Fractional chromatic numbers of planar graphs and $K_5$-minor-free graphs}\label{sec:proofs}

\subsection[odd girth 5]{Fractional colorings of graphs of (odd) girth at least $5$}
We are now ready to apply~\Cref{thm:oct} to prove~\Cref{thm:general}. 
\begin{proof}[Proof of~\Cref*{thm:general}]
Let $G$ be a graph of odd girth at least $\ell$ with $\mad(G)<2\ell$ and assume that $\Poc(G)$ is integral. Let $\boldsymbol{w}\in\R_{\geq 0}^{V(G)}$ be an arbitrary vector. By~\Cref{thm:oct}, there exists an odd-cycle transversal $Z$ such that
$w(Z)\leq\frac{\mad(G)}{2\ell}w(V(G))$. Since $G-Z$ is bipartite, the vertex set of one bipartition of $G-Z$ of larger weight is an independent set $I$ of $G$ satisfying
$$w(I)\geq\frac{w(V(G))-w(Z)}{2}
\geq\frac{2\ell-\mad(G)}{4\ell}w(V(G)).$$
Since $\mad(G)<2\ell$ and $\boldsymbol{w}$ is arbitrary, it follows from~\Cref{lem:criterion}, applied with $r=\frac{4\ell}{2\ell-\mad(G)}$, that $\chi_f(G)\leq\frac{4\ell}{2\ell-\mad(G)}$. This proves~\Cref{thm:general}. 
\end{proof}

The equality holds in~\Cref{thm:general} for odd cycles $C_\ell$ when
$\ell\geq 3$ is odd, since $\mad(C_\ell)=2$ and
$\chi_f(C_\ell)=\frac{2\ell}{\ell-1}$.

\begin{proof}[Proof of~\Cref*{thm:planar}]
Let $k$ be a positive integer and let $G$ be an $n$-vertex planar graph of girth at least $2k+3$, where $n\geq 3$. It is obvious that $G$ is of odd girth at least $2k+3$ and is $K_5$-minor-free. Hence, by \Cref{lem:integral}, $\Poc(G)$ is integral. By the Euler formula, we have that
$\mad(G)\leq \frac{2(2k+3)(n-2)}{(2k+1)n}<\frac{2(2k+3)}{2k+1}<2(2k+3)$. Applying \Cref{thm:general} with $\ell=2k+3$ yields
$$
    \chi_f(G)\leq\frac{4(2k+3)}{2(2k+3)-\mad(G)}\leq 
\frac{4(2k+3)}{2(2k+3)-\frac{2(2k+3)(n-2)}{(2k+1)n}}
    =\frac{(2k+1)n}{kn+1}.
$$
This completes the proof of~\Cref{thm:planar}.
\end{proof}

We now focus on the fractional chromatic number of $K_5$-minor-free graphs. A standard consequence of Wagner's characterization~\cite{Wagner} is
that a $K_5$-minor-free graph $G$ with $n\geq 3$ vertices has at most $3n-6$ edges. Bar\'at~\cite{Barat} proved that a $K_5$-minor-free graph of girth at least $5$ with $n\geq 4$ vertices and $m$ edges satisfies $5m\leq 9n-21$, except for $C_5$ and the Petersen graph with one edge deleted. In these two exceptional graphs, $5m=9n-20$. Therefore, we have the following maximum average degree conditions.

\begin{lemma}\label{lem:minor-density}
Let $n\geq 4$ be a positive integer and let $G$ be an $n$-vertex $K_5$-minor-free graph.
\begin{enumerate}[label=(\arabic*)]
\item{\em \cite{Wagner}}\label{lem:minor-density-six}
$\mad(G)\leq\frac{6n-12}{n}$.
\item{\em \cite{Barat}}\label{lem:minor-density-girth}
If $G$ has girth at least $5$, then $\mad(G)\leq \frac{18n-40}{5n}$.
\end{enumerate}
\end{lemma}

We now prove the corollaries about $K_5$-minor-free graphs.

\begin{proof}[Proof of~\Cref*{thm:K5minorfree}] \Cref{lem:integral} implies that $\Poc(G)$ is integral. By the density hypothesis,
$\mad(G)\leq\frac{2(2k+3)}{2k+1}<2(2k+3)$. \Cref{thm:general}, applied with $\ell=2k+3$, therefore yields
$$\chi_f(G)\leq\frac{4(2k+3)}{2(2k+3)-\mad(G)}\leq \frac{4(2k+3)}{2(2k+3)-\frac{2(2k+3)}{2k+1}}
=\frac{2k+1}{k}.
$$
\end{proof}

\begin{proof}[Proof of~\Cref*{cor:minor-bounds}]
\Cref{lem:integral} implies that $\Poc(G)$ is integral.

\noindent
\ref{cor:minor-odd-girth}
By~\Cref{lem:minor-density}\ref{lem:minor-density-six},
we have $\mad(G)\leq \frac{6n-12}{n}<2\ell$.
\Cref{thm:general} yields
$$\chi_f(G)\leq\frac{4\ell}{2\ell-\mad(G)}\leq\frac{4\ell}{2\ell-\frac{6n-12}{n}}=\frac{2\ell n}{(\ell-3)n+6}.$$

\noindent
\ref{cor:minor-no-four} Since $G$ is of odd girth at least $5$ and has no $4$-cycle, $G$ is of girth at least $5$. 
By \Cref{lem:minor-density}\ref{lem:minor-density-girth}, we have  $\mad(G)\leq\frac{18n-40}{5n}$.
\Cref{thm:general} yields
$$\chi_f(G)\leq\frac{4\ell}{2\ell-\mad(G)}\leq\frac{4\ell}{2\ell-\frac{18n-40}{5n}}=\frac{10\ell n}{(5\ell-9)n+20}.$$
\end{proof}

\subsection{Fractional colorings of triangle-free planar graphs}

We use the following special case of the folding lemma of Klostermeyer and Zhang~\cite{KlostermeyerZhang}; see also Lemma~9 of~\cite{DvorakSereniVolecPlanar}.

\begin{lemma}{\em \cite{KlostermeyerZhang,DvorakSereniVolecPlanar}}
\label{lem:fold}
Let $H$ be a triangle-free plane graph and let $C$ be a facial $4$-cycle. One of the two pairs of opposite vertices of the $4$-cycle $C$ can be identified so that, after parallel edges are suppressed, the resulting graph $H'$ is planar and triangle-free. Moreover, the identification map is a homomorphism of $H$ to $H'$.
\end{lemma}

The next fact is a standard precoloring extension theorem. See, for example, Theorem~5 of~\cite{BorodinKostochkaLidickyYancey}.

\begin{lemma}{\em \cite{BorodinKostochkaLidickyYancey}}
\label{lem:extension}
Let $H$ be a triangle-free plane graph and let $C$ be a facial $4$-cycle. Then each proper $3$-coloring of $C$ can be extended to a proper $3$-coloring of $H$.
\end{lemma}

We also use a basic fact of fractional colorings: If $\chi_f(H)\leq\frac{a}{b}$ for positive integers $a$ and $b$, then there exists a positive integer $M^*$ such that $H$ admits an $(aM\colon bM)$-coloring for every positive integer $M$ divisible by $M^*$.

\begin{proof}[Proof of~\Cref*{thm:triangle-planar}]
Suppose to the contrary that $G$ is a minimum counterexample to \Cref{thm:triangle-planar}. That is to say, $G$ is an $n$-vertex triangle-free plane graph (with a fixed plane embedding) such that $\chi_f(G)>\frac{3n}{n+1}$. By minimality, $G$ is connected. Moreover, it follows from \Cref{thm:planar} that $G$ contains at least one $4$-cycle. 

If $G$ has a facial $4$-cycle, then applying \Cref{lem:fold} we obtain a triangle-free planar graph $G'$ of order $n-1$ and a homomorphism of $G$ to $G'$. By the minimality of $G$, $\chi_f(G')\leq \frac{3(n-1)}{n}$. Hence, $\chi_f(G)\le\chi_f(G')<\frac{3n}{n+1}$, a contradiction. So every $4$-cycle of $G$ is a separating $4$-cycle.

Let $C=abcd$ be a separating $4$-cycle in the fixed plane embedding of $G$. Let $G_1$ and $G_2$ consist of $C$ together with all vertices and edges drawn in its interior and exterior, respectively. Thus,
$G=G_1\cup G_2$ and $G_1\cap G_2=C$. Since $G$ is triangle-free, $C$ has no chord.
Let $n_i:=|V(G_i)|$ for $i\in \{1,2\}$. Then we have $5\leq n_i\leq n-1$ and $n_1+n_2=n+4,$
as each side of $C$ contains at least one vertex.

Note that the cycle $C$ bounds a face in each $G_i$. For each $i\in \{1,2\}$, we can apply
\Cref{lem:fold} to this facial cycle to obtain a triangle-free planar graph $H_i$ of order $n_i-1$ and a homomorphism of $G_i$ to $H_i$. By the minimality of $G$, we have that $\chi_f(G_i)\leq \chi_f(H_i)\leq \frac{3(n_i-1)}{n_i}.$ Let $M$ be a positive integer such that for each $i\in\{1,2\}$, $G_i$ admits
a $(3(n_i-1)M:n_iM)$-coloring $\phi_i$. Define
$$
x_i:=n_iM-|\phi_i(a)\cap\phi_i(c)|,
    \text{~~and~~}
y_i:=n_iM-|\phi_i(b)\cap\phi_i(d)|.
$$
Since every vertex of $\{a,c\}$ is adjacent to every vertex of $\{b,d\}$, the two sets $\phi_i(a)\cup\phi_i(c)$ and $\phi_i(b)\cup\phi_i(d)$ are disjoint. Noting that $|\phi_i(a)\cup\phi_i(c)|=n_iM+x_i$ and $|\phi_i(b)\cup\phi_i(d)|=n_iM+y_i$, $(n_iM+x_i)+(n_iM+y_i)\leq 3(n_i-1)M$. Thus $x_i+y_i\le(n_i-3)M.$ Let $x:=\max\{x_1,x_2\}$ and $y:=\max\{y_1,y_2\}$. For each $i\in\{1,2\}$,
\begin{equation}\label{eq:adjustment}
(x-x_i)+(y-y_i)=\max\{0,x_{3-i}-x_i\}+\max\{0,y_{3-i}-y_i\}\leq x_{3-i}+y_{3-i}\leq (n_{3-i}-3)M.
\end{equation}

To combine the two colorings $\phi_1$ and $\phi_2$, we need them to have the same number of available colors and their restrictions to $C$ to agree after relabeling colors. We therefore enlarge the colorings while making the intersection deficits of the two opposite pairs equal on both sides. We consider proper $3$-colorings on $C$. Let $\psi_0,\psi_1,$ and $\psi_2$ be proper $3$-colorings of $C$ with colors
in $\{\alpha,\beta,\gamma\}$, defined in the vertex order $(a,b,c,d)$ by $\psi_0=(\alpha,\beta,\alpha,\beta),
\psi_1=(\alpha,\beta,\gamma,\beta),$ and $
\psi_2=(\alpha,\beta,\alpha,\gamma).$
By \Cref{lem:extension}, each $\psi_j$ extends to a proper
$3$-coloring $\psi_{i,j}$ of $G_i$ using the same set of colors $\{\alpha,\beta,\gamma\}$.

For each $i\in\{1,2\}$, append to $\phi_i$ a total of $(n_{3-i}-3)M$ such colorings:
$x-x_i$ copies of $\psi_{i,1}$, $y-y_i$ copies of $\psi_{i,2}$, and $(n_{3-i}-3)M-(x-x_i)-(y-y_i)$ copies of $\psi_{i,0}$. These multiplicities are all nonnegative integers by~\eqref{eq:adjustment}. For every copy, replace its three colors by a new three-element set of colors, disjoint from the set of colors used for $\phi_i$ and from the sets used for all other copies.

At each vertex, take the union of its original color set and the colors assigned to it in the appended copies. The resulting coloring $\widehat\phi_i$ is proper. Each vertex receives $n_iM+(n_{3-i}-3)M=(n+1)M$ colors, and the set of available colors has size $3(n_i-1)M+3(n_{3-i}-3)M=3nM$. Thus $\widehat\phi_i$ is a $(3nM:(n+1)M)$-coloring of $G_i$.

The vertices $a$ and $c$ receive the same color in $\psi_0$ and $\psi_2$, but different colors in $\psi_1$. Similarly, $b$ and $d$ receive the same color in $\psi_0$ and $\psi_1$, but different colors in $\psi_2$. Since the sets of colors used for the copies are pairwise disjoint,
$$
|\widehat\phi_i(a)\cap\widehat\phi_i(c)|=n_iM-x_i+(n_{3-i}-3)M-(x-x_i)=(n+1)M-x,$$
$$
|\widehat\phi_i(b)\cap\widehat\phi_i(d)|=n_iM-y_i+(n_{3-i}-3)M-(y-y_i)=(n+1)M-y.
$$
Since adjacent vertices of $C$ receive disjoint color sets, each
available color occurs on $C$ at precisely one of the following
sets of vertices: $\varnothing,\{a\},\{b\},\{c\},
\{d\},\{a,c\},\{b,d\}.$
The classes corresponding to $\{a\}$ and $\{c\}$ each have size $x$, and those corresponding to $\{b\}$ and $\{d\}$ each have size $y$. The two opposite-pair classes have the intersection sizes displayed above, and the class of colors unused on $C$ has size
$(n-2)M-x-y$. This is approperiate since we have $x+y\le x_1+y_1+x_2+y_2\le (n-2)M$ by~\eqref{eq:adjustment}.  Consequently, all seven classes have the same sizes
in $\widehat\phi_1$ and $\widehat\phi_2$.

Choosing a bijection between each pair of corresponding classes gives a relabeling of the colors of $\widehat\phi_2$ under which the two colorings agree at every vertex of $C$. They therefore combine into a $(3nM:(n+1)M)$-coloring of $G$, contradicting the fact that $\chi_f(G)>\frac{3n}{n+1}$.
\end{proof}

\section{Concluding remarks}\label{sec:conclusion}

First we remark that the strict inequality in \Cref{thm:planar} holds for each graph. It is natural to ask whether this can be strengthened to a uniform bound: for every integer $k\geq 1$, does there exist a constant $\varepsilon_k>0$, depending only on $k$, such that every planar graph $G$ of girth at least $2k+3$ satisfies $\chi_f(G)\leq\frac{2k+1}{k}-\varepsilon_k$?
For $k=1$, this is precisely the question posed by
Dvo\v{r}\'ak and Mnich~\cite{DvorakMnich}, also explicitly stated as a conjecture by Dvo\v{r}\'ak, Sereni, and Volec
in~\cite{DvorakSereniVolecPlanar}.

Our original motivation is to resolve \Cref{conj:Naserasr-fractional}. Our proof of its girth version uses the girth assumption to bound the maximum average degree. To extend this approach to the odd-girth setting, one possible direction is to establish stronger weighted bounds for odd-cycle transversals. Specifically, it would suffice to prove that, for every integer $k\geq 1$, every planar graph $G$ of odd girth at least $2k+3$, and every $\boldsymbol{w}\in\R_{\geq 0}^{V(G)}$, there exists an odd-cycle transversal $Z$ satisfying $w(Z)\leq\frac{w(V(G))}{2k+1}.$ 

Within this framework, a natural direction is to exploit the special edge costs $c_{uv}=\min\{w_u,w_v\}$ to obtain sharper odd-cycle edge-cover bounds than those derived from maximum
average degree alone. Another, potentially more flexible, direction is to refine the construction of the independent set. Rather than deleting $Z$ and using only a bipartition class of $G-Z$, one could combine suitably chosen vertices of $Z$ with an independent set in $G-Z$. This could recover some of the
weight lost when deleting the transversal and yield the required independent-set bound without requiring such a strong bound on $w(Z)$.

\section*{Declaration of AI usage}

During the development and preparation of this manuscript, the authors used ChatGPT 6 Sol to explore potential approaches and refine the presentation. All AI-assisted suggestions and arguments were carefully reviewed, revised, and independently verified by the authors. The authors take full responsibility for the mathematical correctness and content of the manuscript.

\section*{Acknowledgments}

Xiaolan Hu is partially supported by National Natural Science Foundation of China (No. 12471325), and the Fundamental Research Funds for the Central Universities (Nos. CCNU26null001, XJ2026000301).
Jiaao Li is partially supported by National Key Research and Development Program of China (No. 2022YFA1006400), National Natural Science Foundation of China (No. 12571371), Natural Science Foundation of Tianjin (No. 24JCJQJC00130), and the Fundamental Research Funds for the Central Universities, Nankai University.
Zhouningxin Wang is partially supported by National Natural Science Foundation of China (Nos. 12301444, 12671413) and the Fundamental Research Funds for the Central Universities, Nankai University.

\end{document}